\documentclass[vecarrow]{svmult}

\usepackage{makeidx}         
\usepackage{graphicx}        
\usepackage{psfrag}                             
\usepackage{multicol}        
\usepackage{amsmath}
\usepackage{amssymb}
\usepackage{amsfonts}
\usepackage{latexsym}
\usepackage[T1]{fontenc}
\usepackage{subfigure}                      
\usepackage{float}
\usepackage{mathrsfs}
\usepackage{mathptmx}         

\begin{document}

\title{Relaxation of Periodic and Nonstandard Growth
	Integrals by means of Two-scale convergence}
\index{Relaxation}

\titlerunning{Two-scale relaxation}

\author{Joel Fotso Tachago, Hubert Nnang and Elvira Zappale}

\authorrunning{J. Fotso Tachago, H. Nnang and E. Zappale}

\institute{J. Fotso Tachago \at 
	University of Bamenda , Faculty of Sciences, Department of Mathematics and Computers Sciences, P.O. Box 39, Bambili, Cameroon\hfill \break
	and \hfill \break
	Dipartimento di Ingegneria Industriale, Universit\'a degli Studi di Salerno, Via Giovanni Paolo II, 132, Fisciano (Sa) 84084, Italy
	 \hfill\break \email{fotsotachago@yahoo.fr},\and
	H. Nnang, \at
	University of Yaound\'e I, and \'Ecole Normale Sup\'erieure de Yaound\'e, P.O. Box 47 Yaounde, Cameroon, \hfill\break \email{hnnang@u1.uninet.cm; hnnang@yahoo.fr}, 
\and 
E. Zappale \at  Dipartimento di Ingegneria Industriale, Universit\'a degli Studi di Salerno, 84084 Fisciano (SA), Italy \hfill\break
\email{ezappale@unisa.it.}}

\maketitle

\index{integral!equations, nonuniqueness}
\index{plane!elasticity}

\bgroup


\renewcommand{\b}{\beta}
\renewcommand{\c}{\gamma}
\newcommand{\n}{\nu}
\newcommand{\x}{\xi}
\renewcommand{\l}{\lambda}
\newcommand{\m}{\mu}
\newcommand{\h}{\chi}
\newcommand{\s}{\sigma}
\renewcommand{\t}{\theta}
\newcommand{\trans}{^{\rm T}}
\newcommand{\ds}{{\partial S}}
\renewcommand{\o}{\omega}
\newcommand{\bb}{\mathcal B}
\newcommand{\rr}{\mathcal R}
\renewcommand{\ll}{\mathcal L}
\renewcommand{\ss}{\mathcal S}
\renewcommand{\gg}{\mathcal G}
\newcommand{\vv}{\mathcal V}
\newcommand{\ww}{\mathcal W}
\newcommand{\Dminus}{(D$^{\omega-}$)}
\newcommand{\Nminus}{(N$^{\omega-}$)}
\newcommand{\Rminus}{(R$^{\omega-}$)}
\newcommand{\Vminus}{{\mathcal V}^{\omega-}}
\renewcommand{\Re}{\mathop\textrm{Re}}
\renewcommand{\Im}{\mathop{\rm Im}}
\newcommand{\mxy}{|x-y|}
\def\Wminus{{\mathcal W}^{\omega-}}
\def\regular#1{C^2(S^{#1})\cap C^1(\bar{S}^{#1})}
\def\holder#1{C^{\,\,#1,\alpha}(\partial S)}
\def\ppd#1{\dfrac{\partial}{\partial#1}}
\newcommand{\e}{\varepsilon}
\newcommand{\supsig}{^{(\sigma)}}
\newcommand{\up}{\upsilon}
\newcommand{\Up}{\Upsilon}
\newcommand{\p}{\psi}
\renewcommand{\a}{\alpha}
\newcommand{\f}{\varphi}
\newcommand{\F}{\Phi}
\renewcommand{\P}{\Psi}
\renewcommand{\d}{\delta}
\renewcommand{\k}{\kappa}
\newcommand{\sums}{\sum_{\s=1}^2}
\newcommand{\summ}{\sum_{m=0}^\infty}
\newcommand{\intl}{\int\limits}
\newcommand{\pd}{\partial}
\renewcommand{\D}{\Delta}

\newcommand{\Cc}{{\mathcal C}}
\newcommand{\Dd}{{\mathcal D}}
\newcommand{\Fill}{\makebox[6pt]{\hfill}}

\newcommand{\tr}{{\mathrm T}}

\allowdisplaybreaks

\abstract{An integral representation result is obtained for the variational limit of the family  functionals $\int_{\Omega}f(\frac{x}{\varepsilon},D u)dx$, $\varepsilon >0$, when the integrand $f = f(x, v)$ is a Carath\'eodory function, periodic in $x$, convex in $v$ and with nonstandard growth.}

\section{Introduction}\label{sec:con1}

In \cite{TN1}, the authors extended the notion of two-scale convergence introduced by \cite{N} (see also \cite{A}, \cite{CDG}, \cite{FZ}, \cite{V} among a wider literature for extensions and related notions) to the Orlicz-Sobolev setting and obtained, under strict convexity assumption on $f$ and suitable boundary conditions, the existence of a unique minimizer for a suitable limit functional as the limit of the minimizers of the original functionals $\int_{\Omega}f(\frac{x}{\varepsilon},D u)dx$ as $\varepsilon \to 0$.

\medskip
In particular they proved (cf. \cite[Corollary 5.2]{TN1}) that for every sequence $(u_\varepsilon)_\varepsilon \in  W^1L^B(\Omega;\mathbb R)$ such that $(Du_\varepsilon)_\varepsilon$ weakly $2s$-converges
to $\mathbb Du_0 = Du +D_yu_1$, where $u_0 = (u,u_1) \in W^{1}L^B(\Omega)\times L^1(\Omega;W^1_\sharp L^B_{\rm per}(Y))$. Then
\begin{equation}
\label{lb2s}\iint_{\Omega \times Y}f(y,\mathbb D u_0)dxdy \leq \liminf_{\varepsilon \to 0}\int_{\Omega}f\left(\frac{x}{\varepsilon},D u_\varepsilon\right)dx,
\end{equation}
where $Y:=(0,1)^d$ ($d \in \mathbb N$) and $\mathbb D u_0:= D u + D_y u_1$ (see Section 2 for the notations adopted in this introduction). 

On the other hand, by the very nature of two-scale converge they obtain, under homogeneous boundary conditions on $\partial \Omega$, (the same proof can be performed for any boundary conditions), for $u$ and $u_1$ regular, the existence of a suitable sequence $(\overline u_\varepsilon)_\varepsilon\subseteq W^1L^B(\Omega)$ such that $\overline u_\varepsilon \rightharpoonup u$, and the opposite inequality holds:
$$
\lim_{\varepsilon\to 0}\int_{\Omega}f\left(\frac{x}{\varepsilon},D \overline u_\varepsilon\right)dx= \iint_{\Omega\times Y}f(y,\mathbb Du_0)dxdy.
$$

Here, by means of two scale convergence we aim to extend their result to any couple of functions $(u_0\equiv)(u, u_1) \in W^{1}L^B(\Omega)\times L^1(\Omega;W^1_\sharp L^B_{\rm per}(Y))$, and also to obtain an integral representation result for
$$
\inf\left\{\liminf_{\varepsilon \to 0}\int_{\Omega}f\left(\frac{x}{\varepsilon}, D u_\varepsilon\right)dx: u_\varepsilon \rightharpoonup u \hbox{ weakly in }W^1L^B(\Omega)\right\}.
$$
Indeed, after stating preliminary results in Section 1.2 on Orlicz-Sobolev spaces and homogenization theory, in Section 1.3 we will prove the following theorem:
\begin{theorem}
\label{THM}
Let $\Omega$ be a bounded open set with Lipschitz boundary and let $f:\Omega \times \mathbb R^d \to \mathbb R$ be a Carath\'eodory function such that
\begin{equation}
\nonumber
f(x,\cdot) \hbox{ is convex for a.e.}x \in \Omega,
\end{equation}
and there exist constants $c, c'$ and $C\in \mathbb R^+$ such that for a.e. $x \in \Omega$  and every $\xi \in \mathbb R^d$, 
\begin{equation}
\label{Bgrowth}
c B'(|\xi|) - c'\leq f(x,\xi)\leq C(1+ B(|\xi|) 
\end{equation}
with $B, B'$ equivalent $N$-functions which satisfy  the $\triangle_2$ condition.
Then, it results that for every $u \in W^{1}L^B(\Omega)$,
\begin{align}
\inf\left\{\liminf_{\varepsilon \to 0}\int_{\Omega}f\left(\frac{x}{\varepsilon},D u_\varepsilon\right)dx: u_\varepsilon \rightharpoonup u \hbox{ weakly in }W^1L^B(\Omega)\right\}
\label{eqthm}\\
=\inf\left\{\limsup_{\varepsilon \to 0}\int_{\Omega}f\left(\frac{x}{\varepsilon},D u_\varepsilon\right)dx: u_\varepsilon \rightharpoonup u \hbox{ weakly in } W^1L^B(\Omega)\right\}
= \int_{\Omega}f_{\rm hom}(Du)dx,\nonumber
\end{align}
where $f_{\rm hom}:\mathbb R^d \to \mathbb R$ is the density defined by
\begin{equation}
\label{fhom}
f_{\rm hom}(\xi):=\inf\left\{\int_{Y}f(y, \xi+D u)dy: u \in W^1 L^B_{\rm per}(Y)\right\}.
\end{equation}
\end{theorem}

We underline that the analysis presented in this paper, holds also in the vectorial case, i.e. fields $u \in W^1L^B(\Omega;\mathbb R^m)$, with the exact same techniques, provided that $f(x,\cdot)$ is convex.

Furthermore, in order to prove \eqref{eqthm}, we also obtain for every $u_0\in W^{1}L^B(\Omega)\times L^1(\Omega;W^1_\sharp L^B_{\rm per}(Y))$, the following two-scale representation:

\begin{align}
\inf\left\{\liminf_{\varepsilon \to 0}\int_{\Omega}f\left(\frac{x}{\varepsilon},D u_\varepsilon\right)dx: u_\varepsilon \overset{2s}{\rightharpoonup} u_0 
\right\}
\nonumber\\
=\inf\left\{\limsup_{\varepsilon \to 0}\int_{\Omega}f\left(\frac{x}{\varepsilon},D u_\varepsilon\right)dx: u_\varepsilon \overset{2s} {\rightharpoonup} u_0 
\right\}\nonumber
= \iint_{\Omega\times Y}f(y, \mathbb Du_0)dxdy.
\end{align}

\section{Preliminaries}
This section is devoted to fix notation adopted in the sequel and state preliminary results on Orlicz-Sobolev spaces and homogenization results that will be exploited in the next section. For more details concerning these latter results, for the sake of brevity, we will refer directly to \cite{TN1}. 

$\Omega\subset \mathbb R^d$ ($d \in \mathbb N$) denotes a bounded open set with Lipschitz boundary. 
\subsection{Orlicz-Sobolev spaces}
Let $B:\left[ 0,+\infty \right[ \rightarrow \left[ 0,+\infty \right[
$ be an $N-$function as in \cite{Ad}, i.e., $B$ is continuous, convex, with $B\left( t\right) >0$ for $t>0,\frac{B\left( t\right) }{t}\rightarrow 0$ as $%
t\rightarrow 0,$ and $\frac{B\left( t\right) }{t}\rightarrow \infty $ as $%
t\rightarrow \infty .$

Equivalently, $B$ is of the form $B\left( t\right)
=\int_{0}^{t}b\left( \tau \right) d\tau ,$ where $b:\left[ 0,+\infty \right[
\rightarrow \left[ 0,+\infty \right[ $ is non decreasing, right continuous,
with $b\left( 0\right) =0,b\left( t\right) >0$ if $t>0$ and $b\left(
t\right) \rightarrow +\infty $ if $t\rightarrow +\infty.$ We denote by $\widetilde{B},$ the Fenchel's conjugate, also called the complementary $N-$function of $B$ defined by $$\widetilde{B}(t)=\sup_{s\geq 0}\left\{ st-B\left( s\right) \right\}, \;\; t\geq 0.
$$ 
It can be proven that (see \cite[Lemma 2.1]{TN1}) if $B$ is an $N$-function and $\tilde B$ is its conjugate, then for all $t>0$, it results
\begin{align}
\frac{t b\left( t\right) }{B\left( t\right) }\geq 1
(> 1 \hbox{ if }b\hbox{ is strictly increasing}),\nonumber\\
\widetilde{B}\left( b\left( t\right) \right) \leq
tb\left( t\right) \leq B\left( 2t\right).
\nonumber
\end{align}%
An $N-$function $B$\ is of class $\triangle _{2}$ (denoted $B\in \triangle
_{2}$) if there are $\alpha >0$ and $t_{0}\geq 0$ such that $B\left(
2t\right) \leq \alpha B\left( t\right) $ for all $t\geq t_{0}$.
 
\medskip 
In all what
follows $\ B$ and $\widetilde{B}$\ are conjuguates $N-$function$s$
satisfying the delta-2 ($\triangle_2$) condition and $c$ refers to a constant that may vary from line to line. 

\noindent The Orlicz-space $L^{B}\left(
\Omega \right) =\left\{ u:\Omega \rightarrow 
\mathbb{C}
\text{ measurable, }\underset{\delta \rightarrow 0^{+}}{\lim }\int_{\Omega
}B\left( \delta \left\vert u\left( x\right) \right\vert \right) dx=0\right\} 
$ is a Banach space for the Luxemburg norm: $$\left\Vert u\right\Vert
_{B,\Omega }=\inf \left\{ k>0:\int_{\Omega }B\left( \frac{\left\vert u\left(
	x\right) \right\vert }{k}\right) dx\leq 1\right\} <+\infty .$$ It follows
that: $C^\infty_c\left( \Omega \right) $ is dense in $L^{B}\left( \Omega
\right) ,L^{B}\left( \Omega \right) $ is separable and reflexive, the dual
of $L^{B}\left( \Omega \right) $ is identified with $L^{\widetilde{B}}\left(
\Omega \right) ,$ and the norm induced  on $L^{\widetilde{B}}\left( \Omega \right) $ as a dual space
is equivalent to $\left\Vert .\right\Vert _{\widetilde{B},\Omega }.$

	
	
 
Analogously one can define the Orlicz-Sobolev functional space as follows: 
$$
W^{1}L^{B}\left( \Omega \right) =\left\{ u\in L^{B}\left( \Omega \right) :%
\frac{\partial u}{\partial x_{i}}\in L^{B}\left( \Omega \right) ,1\leq i\leq
d\right\} ,$$ where derivatives are taken in the distributional sense on $%
\Omega .$ 
Endowed with the norm $\left\Vert u\right\Vert _{W^{1}L^{B}\left(
	\Omega \right) }=\left\Vert u\right\Vert _{B,\Omega }+\sum_{i=1}^{d}$ $%
\left\Vert \frac{\partial u}{\partial x_{i}}\right\Vert _{B,\Omega },u\in
W^{1}L^{B}\left( \Omega \right) ,$\ \ $W^{1}L^{B}\left( \Omega \right) $ is
a reflexive Banach space. We denote by $W_{0}^{1}L^{B}\left( \Omega \right)
, $ the closure of $C^\infty_c\left( \Omega \right) $\ in $%
W^{1}L^{B}\left( \Omega \right) $ and the semi-norm $u\rightarrow \left\Vert
u\right\Vert _{W_{0}^{1}L^{B}\left( \Omega \right) }=\left\Vert
Du\right\Vert _{B,\Omega }=\sum_{i=1}^{d}$ $\left\Vert \frac{\partial u}{%
	\partial x_{i}}\right\Vert _{B,\Omega }$ is a norm on $W_{0}^{1}L^{B}\left(
\Omega \right) $ equivalent to $\left\Vert .\right\Vert _{W^{1}L^{B}\left(
	\Omega \right) }.$

\subsection{Homogenization}
In order to deal with periodic integrands we will adopt the following notation.

Let $Y:=(0,1)^d$. The letter $\varepsilon$ throughout will denote a family of positive real numbers converging to $0$. The set $\mathbb R^d_y$ will denote $\mathbb R^d$, but the subscript $y$ emphasizes the fact that this is the set where the space variable $y$ is.
 We also define
$$C_{\rm per}(Y) =\{v \in C(\mathbb R^d_y): Y-{\rm periodic}\},$$
and
$$L^B_{\rm per}(Y):=\{v \in L^B_{\rm loc}(\mathbb R^N_y): Y-{\rm periodic}\}.
$$
Moreover we observe that $L^B_{\rm per}(Y)$ is a Banach space under the Luxemburg norm 
$\|\cdot\|_{B,Y}$, and $C_{\rm per}(Y)$ is dense in
$L^B_{\rm per}(Y)$ (see \cite[Lemma 2.1]{TN1}). 

For $v \in L^B_{\rm per}(Y)$ 
let
$$v^\varepsilon(x) = v\left(\frac{x}{\varepsilon}\right), 
\hspace{0.05cm} x \in \mathbb R^d.$$
Given $v \in L^B_{\rm loc}(\Omega \times \mathbb R^N_y)$ and $\varepsilon >0$, we put
$$v^\varepsilon (x) = v\left(x, \frac{x}{\varepsilon}\right), 
\hspace{0.05cm} x \in \mathbb R^d \hbox{ whenever it makes sense.}$$
We define the vector space
$$
L^B(\Omega \times Y_{\rm per}):=\{u \in L^B_{\rm loc}(\Omega \times \mathbb R^N_y): \hbox{ for a.e. } x \in \Omega, u(x, \cdot) \hbox{ is } Y-\hbox{ periodic}\}.
$$
and observe that the embbedding $L^B(\Omega, C_{\rm per}(Y)) \to L^B(\Omega \times Y_{\rm per})$
is continuous.

Moreover we will make use of the space
\begin{align}
\nonumber
W^{1}L^B_{\rm per}(Y):=\left\{ u \in W^{1}L^B_{\rm loc}(\mathbb R^N_y): u, \, \frac{\partial u}{\partial x_i}, i=1,\dots, N, \;  Y-\hbox{ periodic}\right\}
\end{align}
where the derivative $\frac{\partial u}{\partial x_i}$
is taken in the distributional sense on $\mathbb R^N_y$, and we endow it 
with the norm
$\|u\|_{W^1L^B_{\rm per}}= \|u\|_{B,Y}+ \sum_{i=1}^N\left\|\frac{\partial u}{\partial x_i}\right\|_{B,Y}$, which makes it a Banach space.

We also consider the space
$$W^1_\sharp L^B_{\rm per}(Y)=\left\{u \in W^1L^B_{\rm per}(Y) : \int_{Y}u(y)dy  =0\right\},$$
and we endow it with the gradient norm

$$\|u\|_{W^1_\sharp L^B_{\rm per}(Y)}= \sum_{i=1}^N\left\|\frac{\partial u}{\partial x_i}\right\|_{B,Y}.$$
Denoting by $C^\infty_{\rm per}(Y)= C_{\rm per}(Y)\cap  C^\infty(\mathbb R^N)$, and recalling that the space
$
C^\infty_{\sharp,{\rm per}}(Y ;\mathbb R) = \left\{u \in C^\infty_{\rm per}(Y ;\mathbb R) \int_Y: u(y)dy = 0\right\}$ is dense in $W^1_\sharp L^B_{\rm per}(Y)
$, one can deduce (cf. \cite{TN1}) the density of the embedding
\begin{equation}
\label{densityTN}
C^\infty_c(\Omega;\mathbb R) \otimes C^\infty_{\sharp,{\rm per}}(Y ;\mathbb R)\subseteq L^1(\Omega;W^1_\sharp L^B_{\rm per}(Y)).
\end{equation}

\bigskip In \cite{TN1} the notion of two-scale convergence introduced by \cite{N} and developed by \cite{A} (see also, among a wide literature, \cite{CDG}, \cite{FZ}, \cite{Ne}, \cite{V} for further developments and related notions like periodic unfolding method), has been extended to the Orlicz-Sobolev setting.

\begin{definition}
	A sequence of functions $\left( u_{\varepsilon }\right) _{\varepsilon}$ in $L^{B}\left( \Omega \right) $ is said to be:
	
	-weakly two-scale convergent in $L^{B}\left( \Omega \right) $ to a function $u_{0}\in L^{B}\left( \Omega \times Y_{\rm per}\right) $ if for every $ \varepsilon \rightarrow 0,$ we have 
	\begin{equation}\label{2bis}
	\int_{\Omega }u_{\varepsilon }f^{\varepsilon }dx\rightarrow \iint_{\Omega \times Y}u_{0}fdxdy,\text{ for all }f\in L^{%
		\widetilde{B}}\left( \Omega ;C_{per}\left( Y\right) \right)
	\end{equation}%
	
	- stongly two-scale convergent in $L^{B}\left( \Omega \right) $\
	to $u_{0}\in L^B\left( \Omega \times Y_{\rm per}\right) $\ if for $%
	\eta >0$ and $f\in L^{B}\left( \Omega ;C_{per}(Y) \right) $ verifying $\left\Vert u_{0}-f\right\Vert _{L^B(\Omega
		\times Y)}\leq \frac{\eta }{2}$ there exist $\rho >0$ such that $%
	\left\Vert u_{\varepsilon }-f^{\varepsilon }\right\Vert _{L^B(\Omega) }\leq
	\eta $ for all $0<\varepsilon \leq \rho .$
\end{definition}
When \eqref{2bis} happens for all $f\in L^{\widetilde{B}}\left( \Omega ;C_{per}\left( Y\right) \right) $ we denote it by "$%
u_{\varepsilon }\overset{2s}{\rightharpoonup} u_{0}$ in $L^{B}\left( \Omega
\right) $ weakly" or simply " $u_{\varepsilon }\overset{2s}{\rightarrow} u_{0}$ in $L^{B}\left( \Omega \right) $\ two-scale weakly" and we will say that 
$u_{0}$\ is the weak two-scale limit in $L^{B}\left( \Omega
\right) $\ \ of the sequence $\left( u_{\varepsilon }\right) _{\varepsilon}.$
In order to denote strong two scale convergence of $u_\varepsilon \to u_0$ we adopt the symbol $\|u_\varepsilon - u_0\|_{2s-L^B(\Omega\times Y)}\to 0$.

The following result, whose proof can be found in \cite{TN1}, allows to extend the notion of weak two-scale convergence at Orlicz-Sobolev functions, guaranteeing, at the same time, a compactness result.

\begin{proposition}
	Let $\Omega $ be a bounded open set in $\mathbb{R}^d$ and let $\left( u_{\varepsilon
	}\right) _{\varepsilon}$ be bounded in $W^{1}L^{B}\left( \Omega \right).$ 
	There exist a subsequence, still denoted in the same way, and $u\in W^{1}L^{B}\left(
	\Omega \right) ,$ $ u_{1}\in L^{1}\left(
	\Omega ;W_{\#}^{1}L^{B}_{\rm per}\left( Y\right) \right)$ such that:
	
	$\left( i\right) $ $u_{\varepsilon }\overset{2s}{\rightharpoonup} u$ in $L^{B}\left(
	\Omega \right) $,
	
	$\left( ii\right) $ $D_{x_{i}}u_{\varepsilon }\overset{2s}{\rightharpoonup}
	D_{x_{i}}u+D_{y_{i}}u_1$ in $L^{B}\left( \Omega
	\right) $, $1\leq i\leq d$.
\end{proposition}

In the sequel we denote by $u_0(x,y)$ the function $u(x)+ u_1(x,y)$, and by $\mathbb D u_0$ the vector
$D u+ D_y u_1$.

\noindent For the sake of brevity, we cannot explicitly quote all the results used throughout the paper but we will refer to \cite{TN1} for futher necessary properties of Orlicz-Sobolev spaces, two-scale convergence and homogenization in the Orlicz setting.

\section{Proof of Theorem \ref{THM}}\label{sec}

This section is devoted to the proof of Theorem \ref{THM}. To this aim recall the definition of $f_{\rm hom}$ given by \eqref{fhom}. 

\begin{proof}[of Theorem \ref{THM}]
	We start observing that the coercivity assumptions on $f$, the compactness result, given by Proposition 1, and \eqref{lb2s} guarantee that 
 for every $u_\varepsilon \rightharpoonup u\in W^{1}L^B(\Omega)$
 \begin{align}\label{lb2ssec}
 \liminf_{\varepsilon \to 0}\int_{\Omega}f\left(\frac{x}{\varepsilon},D u_\varepsilon\right)dx
 \geq \iint_{\Omega \times Y}f(y,\mathbb D u_0)dxdy, 
 \end{align}
 where $u_0(x,y) = u(x)+u_1(x,y)$ is the weak two scale limit of $(u_\varepsilon)_\varepsilon$.
 Clearly passing to the infimum on both sides of the equation \eqref{lb2ssec}, and recalling that 
 $
 \mathbb D u_0= D u+ D_y u_1,
 $ we obtain
  \begin{align}\nonumber
 \inf\left\{\liminf_{\varepsilon \to 0}\int_{\Omega}f\left(\frac{x}{\varepsilon},D u_\varepsilon\right)dx: u_\varepsilon \rightharpoonup u \hbox{ in }W^1L^B(\Omega)\right\}
 \\
 \geq\inf\left\{\iint_{\Omega \times Y}f(y, Du+ D_y u_1)dxdy: u_1 \in L^1(\Omega;W^1_\sharp L^B_{\rm per}(Y))\right\}=\int_\Omega f_{\rm hom}(Du)dx,\nonumber
 \end{align}
 where one can replicate the same proof as \cite[Lemma 2.2]{CDDA} replacing $t$ by $1$ and $f_1$ by $f_{\rm hom}$ in \eqref{fhom} therein and exploit the convexity of $f$ to replace functions with null boundary datum on $\partial Y$, with periodic ones (see also end of \cite[Chapter 3]{Ne}).

\noindent The upper bound exploits an argument very similar to the one presented in \cite{Ne}, relying, in the present context, on the density result in \eqref{densityTN}. Indeed we can first observe that, as in \cite[Corollary 5.1]{TN1} 
for any given $u\in  C^\infty(\overline \Omega)$ and $\phi_1 \in C^\infty_c(\Omega)\otimes C^\infty_{\rm per}(Y)$ it results that, given $\phi_\varepsilon(x):=u(x)+ \varepsilon \phi_1\left(x, \frac{x}{\varepsilon}\right)$
\begin{equation}
\label{ub2s1}
\lim_{\varepsilon\to 0}\int_{\Omega}f\left(\frac{x}{\varepsilon}, D\phi_\varepsilon(x)\right)dx = \iint_{\Omega \times Y}f(y,Du +D_y\phi_1(x,y))dxdy.
\end{equation}
On the other hand, given $u \in W^1L^B(\Omega)$ and $u_1 \in L^1(\Omega;W^1_\sharp L^B_{\rm per}(Y))$, \eqref{densityTN} guarantees that for each $\delta > 0$ we can find maps
$u_\delta \in C^\infty(\overline \Omega)$ and $v_\delta \in
 C^\infty_c(\Omega;C^\infty_{\rm per}(Y))$ (this latter with zero averageb)
such that
\begin{equation}
\label{Bdense}
\|u-u_\delta\|_{W^1L^B(\Omega)}+\|u_1-v_\delta\|_{L^1(\Omega;W^{1}L^B_{\rm per}(Y)}\leq \delta
\end{equation}

Next defining, for every $\delta$, and for every $x \in \Omega$, 
$$
u_{\delta,\varepsilon}(x):= u_\delta(x)+ \varepsilon v_\delta^\varepsilon(x),
$$
one has
$$
Du_{\delta, \varepsilon}(x)=Du_\delta(x)+ \varepsilon D_x v_\delta\left(x, \frac{x}{\varepsilon}\right)+ D_y v_\delta\left(x, \frac{x}{\varepsilon}\right).
$$
Clearly, as $\varepsilon \to 0$, it results 
\begin{align*}
u_{\delta, \varepsilon}\rightarrow u_\delta \hbox{ in }L^B(\Omega),\\
Du_{\delta, \varepsilon}(x)\overset{2s}{\rightarrow} Du_\delta(x)+ D_yv_\delta(x,y) \hbox{ strongly in }  L^B(\Omega \times Y_{\rm per}). 
\end{align*}

Now we define
\begin{equation}
\label{cdeltavarepsilon}
c_{\delta,\varepsilon}:=\|u_{\delta, \varepsilon}-u\|_{W^1L^B(\Omega)}+
\left|\left\|D u_{\delta, \varepsilon}\right\|_{L^B(\Omega)}- \left\| Du - D_y u_1\right\|_{L^B(\Omega \times Y)}\right|,
\end{equation}
having the aim of constructing, via a diagonalizing argument, a sequence stongly two scale convergent to $u_0=u+ u_1$.

Thus, it is easily seen that
\begin{equation}
\lim_{\delta\to 0}\lim_{\varepsilon\to 0}c_{\delta, \varepsilon}=0,\nonumber
\end{equation}
which allows us to apply H. Attouch Diagonalization Lemma, thus detecting a sequence $\delta(\varepsilon)\to 0$ as $\varepsilon \to 0$, such that $c_{\delta(\varepsilon),\varepsilon}\to 0$ and $u_{\delta(\varepsilon),\varepsilon}\to u $ in $L^B(\Omega)$, with
$$
D u_{\delta(\varepsilon),\varepsilon}(x)\overset{2s}{\to} Du(x) + D_y u_1(x,y)  \hbox{ strongly in }L^B(\Omega \times Y_{\rm per}).
$$
This latter convergence, and \cite[Remark 4.1]{TN1} ensure that $Du_{\delta(\varepsilon),\varepsilon}\rightharpoonup Du$ weakly
in $L^B(\Omega)$, thus, \eqref{ub2s1}, the continuity of $f$ in the second variable, \eqref{Bgrowth}, guarantee that for every $u_1 \in L^1(\Omega;W^{1}_\sharp L^B_{\rm per}(Y))$,
\begin{align*}
\lim_{\varepsilon\to 0}\int_{\Omega}f\left(\frac{x}{\varepsilon}, Du_{\delta(\varepsilon),\varepsilon}\right)dx = \iint_{\Omega \times Y}f(y,Du +D_yu_1(x,y))dxdy.
\end{align*}
as desired.

 Thus we can conclude that 
 \begin{align}\nonumber
 \inf\left\{\limsup_{\varepsilon \to 0}\int_{\Omega}f\left(\frac{x}{\varepsilon}, D u_\varepsilon\right)dx: u_\varepsilon \rightharpoonup u \hbox{ in }W^1L^B(\Omega)\right\}
 \\
 \leq \lim_{\varepsilon\to 0}\int_{\Omega}f\left(\frac{x}{\varepsilon}, D u_{\delta(\varepsilon),\varepsilon}\right)dx \leq \iint_{\Omega \times Y}f(y, Du(x)+ D_y u_1(x,y)dx dy. \nonumber
 \end{align}  
 Hence
 \begin{align*}
 \inf\left\{\limsup_{\varepsilon \to 0}\int_{\Omega}f\left(\frac{x}{\varepsilon}, D u_\varepsilon\right)dx: u_\varepsilon \rightharpoonup u \hbox{ in }W^{1,B}(\Omega;\mathbb R)\right\}
 \\
 \leq \inf\left\{\iint_{\Omega \times Y}f(y, Du+ D_y u_1)dxdy: u_1 \in L^1(\Omega;W^1_\sharp L^B_{\rm per}(Y))\right\}.
 \end{align*}
 This, together with the last equality in \eqref{lb2ssec} concludes the proof.
\end{proof}

{\bf Acknowledgements}
The first and the third author aknowledge the support of the Programme ICTP-INdAM research in pairs 2018.
Joel Fotso Tachago thanks Dipartimento di Ingegneria Industriale at University of Salerno for its hospitality. 
Elvira Zappale is a member of GNAMPA-INdAM, whose support is gratefully acknowledged.

\end{document}